\documentstyle[12pt]{amsart}
\begin{document}
\title{Biharmonic maps into symmetric spaces and integrable systems}
\author{Hajime URAKAWA}
\address{Division of Mathematics\\
Graduate School of Information Sciences\\
Tohoku University\\
Aoba 6-3-09, Sendai, 980-8579, Japan}
\curraddr{Institute for International Education\\ Tohoku University\\
Kawauchi 41, Sendai, 980-8576, Japan}
\title[Biharmonic maps into symmetric spaces]
{Biharmonic maps into symmetric spaces and integrable systems}
\email{urakawa@@math.is.tohoku.ac.jp}
\keywords{harmonic map, biharmonic map, symmetric space, integrable system, Maurer-Cartan form}
\subjclass[2000] 
{58E20}
\thanks{Supported by the Grant-in-Aid for the Scientific Research, 
(C), No. 21540207, 
Japan Society for the Promotion of Science.}
\dedicatory{}
\maketitle
\begin{abstract} 
In this paper, the description 
of biharmonic map equation 
in terms of the Maurer-Cartan form 
for all smooth map 
of a compact Riemannian manifold into a Riemannian 
symmetric space $(G/K,h)$ induced from the bi-invariant Riemannian metric $h$ on $G$ is obtained. By this formula, 
all biharmonic curves into symmetric spaces are determined, and 
all the biharmonic maps 
of an open domain of ${\Bbb R}^2$ with the standard Riemannian metric into $(G/K,h)$ are characterized. 
\end{abstract}
\numberwithin{equation}{section}
\theoremstyle{plain}
\newtheorem{df}{Definition}[section]
\newtheorem{th}[df]{Theorem}
\newtheorem{prop}[df]{Proposition}
\newtheorem{lem}[df]{Lemma}
\newtheorem{cor}[df]{Corollary}
\newtheorem{rem}[df]{Remark}
\section{Introduction and statement of results.} 
This paper is a continuation of our previous one \cite{U}. 
In our previous paper, we discussed about the description of biharmonic maps into compact Lie groups in terms of the Maurer-Cartan form, and gave their explicit constructions. In this paper, we want to extend them 
to biharmonic maps into Riemannian symmetric spaces. 
\par
The theory of harmonic maps into Lie groups, symmetric spaces 
or homogeneous spaces has been extensively studied 
related to the integrable systems by many authors
(for instance, 
 \cite{Uh}, \cite{W}, \cite{DPW}, \cite{W3}, 
 \cite{MHO1}, \cite{MHO2}, \cite{MO}, \cite{O}, \cite{DSU}). 
 In particular, the moduli space of harmonic maps of 2-sphere 
 into symmetric spaces was completely determined (cf. \cite{DPW}, 
 \cite{GO}, \cite{O}). 
 Let us recall the loop group formulation of  harmonic maps into symmetric spaces, briefly. 
 Let $\varphi$ be a smooth map of a Riemann surface 
 $M$ into a Riemannian symmetric space $(G/K, h)$ with a lift $\psi:\,M\rightarrow G$ 
 so that $\pi\,\circ\, \psi=\varphi$. 
 Let
 $\frak g=\frak k \oplus\frak m$ be the corresponding Cartan decomposition of the Lie algebra $\frak g$ of the Lie group $G$. 
 Then, the pull back $\alpha=\psi^{-1}d\psi$ of the Maurel-Cartan form on $G$ is decomposed as $\alpha=\alpha_{\frak k}+\alpha_{\frak m}$, correspondingly. 
 Let us decompose $\alpha_{\frak m}$ into 
 the sum of holomorphic part and the anti-holomorphic one: 
 $\alpha_{\frak m}=\alpha_{\frak m}{}'+\alpha_{\frak m}{}''$. 
 Then, 
 one can obtain the extended solution $\widetilde{\psi}$ of $M$ into 
 a loop group $\Lambda G$ satisfying that 
 $\widetilde{\psi}^{-1}d\widetilde{\psi}
 =\lambda\alpha_{\frak m}{}'+\alpha_{\frak k}+\lambda^{-1}\alpha_{\frak m}{}''$ 
 for all $\lambda\in U(1)=\{\lambda\in {\mathbb C}:\,\vert\lambda\vert=1\}$. 
 (cf. \cite{DPW}). Then, $\varphi:\, M\rightarrow (G/K,h)$ 
 is harmonic if and only if 
 there exists a holomorphic and horizontal map $\widetilde{\psi}$ 
 of $M$ into the homogeneous $\Lambda G/K$ with 
 $\widetilde{\psi}_1=\psi$ (cf. \cite{DPW}, p. 648).  Then, one can obtain 
 Weierstrass-type representation of harmonic maps (cf. \cite{DPW}, pp. 648--662). 
 \par
 On the other hand, the notion of harmonic map has been extended 
 to the one of biharmonic map (cf. \cite{EL}, \cite{J}).  
 In this paper, we wil  
 describe
 biharmonic maps into Riemannian symmetric spaces 
 in terms of the pull back $\alpha=\alpha_{\frak k}+\alpha_{\frak m}$ of the Maurer-Cartan form 
 (cf. Theorem 3.6),  
 give some explicit solutions of the biharmonic map equation in Riemannian symmetric spaces, and construct several biharmonic maps into Riemannian symmetric spaces (Sections 4 and 5). 
\vskip0.6cm\par
{\bf Acknowledgement}: 
The author expresses his gratitude to Prof. J. Inoguchi and Prof. Y. Ohnita 
who gave many useful suggestions and discussions, and 
Prof. A. Kasue for his financial support 
during the preparation of this paper.  
\vskip0.6cm\par
\section{Preliminaries.}
In this section, we prepare general 
materials and  facts 
harmonic maps, biharmonic maps into Riemannian symmetric spaces (cf. \cite{KN}). 
\subsection
\quad 
Let $(M,g)$ be an $m$-dimensional compact Riemannian manifold, and 
the target space $(N,h)$ is 
an $n$-dimensional Riemannian symmetric space $(G/K,h)$. 
Nemely, 
let 
$\frak g$, $\frak k$ be the Lie algebras of $G$, $K$, and  
${\frak g}={\frak k}\oplus {\frak m}$ is the Cartan decomposition of $\frak g$, 
and $h$, 
the $G$-invariant Riemannian metric on $G/K$ 
corresponding to the Ad$(K)$-invariant inner product 
$\langle\,,\,\rangle$ on $\frak m$. 
Let $k$ be a left invariant Riemannian metric on $G$ such as
the natural projection $\pi:\,G\rightarrow G/K$ is a Riemannian submersion 
of $(G,k)$ onto $(G/K,h)$. 
For every $C^{\infty}$ map $\varphi$ of 
$M$ into $G/K$, let us take its (local) {\em lift} $\psi: M\rightarrow G$ 
of $\varphi$, i.e., 
$\varphi=\pi\,\circ\,\psi$, 
$\varphi(x)=\psi(x)\,K\in G/K$ 
$(x\in U\subset M)$, where $U$ is an open subset of $M$.  
\par
The {\it energy functional} 
on the space $C^{\infty}(M,G/K)$ of all $C^{\infty}$ maps 
of $M$ into $G/K$ is defined by 
$$E(\varphi)=\frac12\int_M\vert d\varphi\vert^2\,v_g,$$
and 
for a $C^{\infty}$ one parameter deformation 
$\varphi_t\in C^{\infty}(M,G/K)$ $(-\epsilon<t<\epsilon)$ 
of  $\varphi$ with $\varphi_0=\varphi$, 
the {\it first variation formula} is given by 
$$
\frac{d}{dt}
\bigg\vert _{t=0}E(\varphi_t)
=-\int_M
\langle 
\tau(\varphi),V \rangle \,v_g,
$$ 
where $V$ is a variation vector field along $\varphi$ defined by 
$V=\frac{d}{dt}\big\vert_{t=0} 
\varphi_t$ 
which belongs to the space 
$\Gamma(\varphi^{-1}T(G/K))$ of sections of the induced bundle of the tangent bundle $T(G/K)$ by $\varphi$. 
The {\it tension field} $\tau(\varphi)$ is 
defined by 
\begin{equation}
\tau(\varphi)=\sum_{i=1}^mB(\varphi)(e_i,e_i),
\end{equation}
where
$$B(\varphi)(X,Y)=\nabla^h_{d\varphi(X)}d\varphi(Y)
-d\varphi(\nabla_XY)$$
for $X,Y\in {\frak X}(M)$. 
Here, $\nabla$, and $\nabla^h$, are the Levi-Civita connections of $(M,g)$ and $(G/K,h)$, respectively. 
For a harmonic map $\varphi:\,(M,g)\rightarrow (G/K,h)$, 
the {\it second variation formula} of the energy functional 
$E(\varphi)$ is 
$$
\frac{d^2}{dt^2}\bigg\vert_{t=0}
E(\varphi_t)=\int_M\langle J(V),V\rangle\,v_g
$$
where 
\begin{align}
J(V)&:=\overline{\Delta}V-{\mathcal R}(V),\\
\overline{\Delta}V&:=\overline{\overline{\nabla}}^{\ast}\,\overline{\overline{\nabla}}V
=
-\sum_{i=1}^m
\{
\overline{\overline{\nabla}}_{e_i}(\overline{\overline{\nabla}}_{e_i}V)
-\overline{\overline{\nabla}}_{\nabla_{e_i}e_i}V
\},\\
{\mathcal R}(V)&:=\sum_{i=1}^m
R^h(V,d\varphi(e_i))d\varphi(e_i).
\end{align}
Here, $\overline{\overline{\nabla}}$ is the induced connection 
on the induced bundle $\varphi^{-1}T(G/K)$, and 
is 
$R^h$ is the curvature tensor of $(G/K,h)$ 
given by 
$R^h(U,V)W=
[\nabla^h_U,\nabla^h_V]W-\nabla^h_{[U,V]}W$ 
$(U,V,W\in {\frak X}(G/K)$).  
\par
The {\it bienergy functional} is 
defined by 
\begin{equation}
E_2(\varphi)=
\frac12
\int_M\vert (d+\delta)^2\varphi\vert^2\,v_g
=\frac12\int_M\vert\tau(\varphi)\vert^2\,v_g, 
\end{equation}
and the {\it first variation formula} of the bienergy 
is given (cf. \cite{J}) by 
\begin{equation}
\frac{d}{dt}\bigg\vert _{t=0}E_2(\varphi_t)
=-\int_M\langle \tau_2(\varphi),V\rangle\,v_g
\end{equation}
where the {\it bitension field} $\tau_2(\varphi)$ is 
defined by 
\begin{equation}
\tau_2(\varphi)=J(\tau(\varphi))
=\overline{\Delta}\tau(\varphi)-{\mathcal R}(\tau(\varphi)),
\end{equation}
and a $C^{\infty}$ map
$\varphi:(M,g)\rightarrow (G/K,h)$ is called to be 
{\it biharmonic} if 
\begin{equation}
\tau_2(\varphi)=0.
\end{equation}
\vskip0.6cm\par
\subsection
\quad 
Let $k$ be a left invariant Riemannian metric on $G$ corresponding to
the inner product $\langle\cdot,\cdot\rangle$ on $\frak g$ given by
$\langle \cdot,\cdot\rangle=-B(\cdot,\cdot)$ 
if $(G/K, h)$ is of compact type, and by
$\langle U+X,V+Y\rangle=-B(U,V)+B(X,Y)$
$(U,V\in {\frak k},\,X,Y\in {\frak m})$ if $(G/K,h)$ is of non-compact type. 
Here, $B(\cdot,\cdot)$ is the Killing form of $\frak g$.   
Then, the projection $\pi$ 
of $G$ onto $G/K$ is a Riemannian submersion of 
$(G,k)$ onto $(G/K,h)$,  
and we have also the orthogonal decomposition of 
the tangent space 
$T_{\psi(x)}G$ $(x\in M)$ with respect to 
the inner product $k_{\psi(x)}(\cdot,\cdot)$ $(x\in M)$ 
in such a way that 
\begin{equation}
T_{\psi(x)}G=V_{\psi(x)}\oplus H_{\psi(x)},
\end{equation}
where the {\em vertical space} at $\psi(x)\in G$ is given by 
\begin{equation}
V_{\psi(x)}={\rm Ker}(\pi_{\ast\,\psi(x)})=\{X_{\psi(x)}\vert\,X\in {\frak k}\},
\end{equation}
and the {\em horizontal space} at $\psi(x)$ is given by 
\begin{equation}
H_{\psi(x)}=\{Y_{\psi(x)}\vert\,Y\in {\frak m}\},
\end{equation}
corresponding to the Cartan decomposition  
${\frak g}={\frak k}\oplus {\frak m}$. 
Then, for every 
$C^{\infty}$ section $W\in \Gamma(\psi^{-1}TG)$, 
we have the decomposition corresponding to (2.9), 
\begin{equation}
W(x)=W^V(x)+W^H(x)\quad (x\in M),
\end{equation}
where $W^V$, $W^H$, (denoted also by 
${\mathcal V} W$, 
${\mathcal H} W$, respectively) belong to $\Gamma(\psi^{-1}TG)$. 
We denote by 
$\Gamma(E)$, the space of all $C^{\infty}$ sections of a vector bundle $E$. 
For $Y\in {\frak m}$, define 
$\widetilde{Y}\in \Gamma(\psi^{-1}TG)$ by 
$
\widetilde{Y}(x):=Y_{\psi(x)}
\quad (x\in M).
$
Let 
$\{X_i\}_{i=1}^n$ be an orthonormal basis of $\frak m$ 
with respect to the inner product $\langle\cdot,\cdot\rangle$ of $\frak g$ corresponding to the 
left invariant Riemannian metric $k$ on $G$. 
Then, $W^H$ can be written in terms of $\widetilde{X_i}$ as
$$
W^H=\sum_{i=1}^nf_i\,\widetilde{X_i}
$$
where $f_i\in C^{\infty}(M)$ $(i=1,\cdots,n)$. 
Because, for every $x\in M$, $W^H(x)\in H_{\psi(x)}$, so that we have
$$
W^H(x)=\sum _{i=1}^n
f_i(x)\,X_{i\,\psi(x)}=\sum_{i=1}^nf_i(x)\widetilde{X_i}(x). 
$$
\par
We say 
$W\in \Gamma(\psi^{-1}TG)$ and $V\in \Gamma(\varphi^{-1}T(G/K))$ 
are {\em $\pi$-related}, denoted by $V=\pi_{\ast}W$, if it holds that 
$$
V(x)=\pi_{\ast}W(x)\quad (x\in M),
$$
where 
$
\pi_{\ast}:\,T_{\psi(x)}G\rightarrow T_{\varphi(x)}(G/K)=T_{\pi(\psi(x))}(G/K)
$
is the differentiation of the projection $\pi$ of $G$ onto $G/K$ at $\psi(x)$ for each $x\in M$. 
\par
Let be $\nabla$, $\nabla^k$, $\nabla^h$, 
the Levi-Civita connections of 
$(M,g)$, $(G,k)$, $(G/K,h)$, and 
$\overline{\nabla}$, 
$\overline{\overline{\nabla}}$, 
the induced connection of $\nabla^k$ on the induced bundle 
$\psi^{-1}TG$ by $\psi:\,M\rightarrow G$, and 
the one of $\nabla^h$ on the induced bundle $\varphi^{-1}T(G/K)$ 
by $\varphi:\,M\rightarrow G/K$, respectively. 
\begin{lem}
Assume that $W\in \Gamma(\psi^{-1}TG)$ and 
$V\in \Gamma(\varphi^{-1}T(G/K))$ are $\pi$-related, i.e., $V=\pi_{\ast}W$. 
\par
$(1)$ Then, we have
\begin{equation}
\overline{\overline{\nabla}}_XV
=\pi_{\ast}
\nabla^k_{(\psi_{\ast}X)^H}W^H,
\end{equation}
where 
$(\psi_{\ast}X)^H$ is the horizontal component of 
$\psi_{\ast}X$ for every $C^{\infty}$ 
vector field $X$ on $M$. 
\par
$(2)$ If we express 
$
W^H=\sum_{i=1}^nf_i\,\widetilde{X_i}
$
and 
$(\psi_{\ast}X)^H=\sum_{j=1}^ng_j\,\widetilde{X_j}$ where 
$f_i, \,g_j\in C^{\infty}(M)\,(i,j=1,\cdots,n)$, 
then, it holds that 
\begin{align}
\left(
\nabla^k_{(\psi_{\ast}X)^H}W^H
\right)_{\psi(x)}
&=
\frac12\sum_{i,j=1}^nf_i(x)\,g_j(x)\,[X_j,X_i]_{\psi(x)}
+
\sum_{i=1}^n
X_x(f_i)\,\widetilde{X_i}(x)\nonumber\\
&\in V_{\psi(x)}\oplus H_{\psi(x)}
\quad(x\in M),
\end{align}
correspondingly. 
\par
$(3)$ For every $x\in M$, we have
\begin{equation}
\overline{\overline{\nabla}}_XV(x)
=\sum_{i=1}^n
X_x(\langle W,X_{i\,\psi(x)}\rangle)\,\pi_{\ast}(X_{i\,\psi(x)}). 
\end{equation}
Here, it holds that 
$\pi_{\ast}(X_{\psi(x)})=t_{\psi(x)\,\ast}\pi_{\ast}(X)$ ($X\in {\frak m}$), 
where $t_a$ is the translation of $G/K$ 
by $a\in G$, i.e., $t_a(yK):=ayK$ $(y\in G)$. 
\end{lem}
\vskip0.6cm\par
\begin{pf}
(1) Due to 
Lemmas 1 and 3 in \cite{O}, p.460, we have 
\begin{align*}
\overline{\overline {\nabla}}_XV&=
\nabla^h_{\varphi_{\ast}X}V\\
&=\nabla^h_{\pi_{\ast}(\psi_{\ast}X)}\pi_{\ast}W\\
&=\pi_{\ast}\left(
{\mathcal H}\,\nabla^k_{(\psi_{\ast}X)^H}W^H\right)\\
&=\pi_{\ast}\nabla^k_{(\psi_{\ast}X)^H}W^H. 
\end{align*}
(2) Indeed, we have 
\begin{align*}
\left(\nabla^k_{(\psi_{\ast}X)^H}W^H
\right)_{\psi(x)}
&=\sum_{j=1}^ng_j(x)\left(\nabla^k_{X_j}W^H
\right)_{\psi(x)}\\
&=\sum_{j=1}^ng_j(x)\left(
\sum_{i=1}^n\nabla^k_{X_j}(f_i\,\widetilde{X_i})
\right)_{\psi(x)}\\
&=
\sum_{i,j=1}^n
g_j(x)
\left\{
(X_jf_i)(x)\,\widetilde{X_i}(x)+f_i(x)\left(
\nabla^k_{X_j}X_i
\right)_{\psi(x)}
\right\}\\
&
=\sum_{i,j=1}g_j(x)
\left\{
(X_jf_i)(x)\,\widetilde{X_i}(x)+\frac12f_i(x)[X_j,X_i]_{\psi(x)}
\right\}\\
&=\frac12\sum_{i,j=1}^nf_i(x)g_j(x)[X_j,X_i]_{\psi(x)}
+\sum_{i=1}^n
X_x(f_i)\,\widetilde{X_i}(x),
\end{align*}
since it holds that  
\begin{align*}
\left(\nabla^k_{X_j}X_i
\right)_{\psi(x)}&=L_{\psi(x)\,\ast}\left(\nabla^k_{X_j}X_i
\right)_e\\
&=L_{\psi(x)\,\ast}\left(
\frac12\,[X_j,X_i]_e
\right)\\
&=\frac12[X_j,X_i]_{\psi(x)}
\end{align*}
and $[{\frak m},{\frak m}]\subset{\frak k}$.  
For (3), notice that 
$W^H=\sum _{i=1}^n\langle W,X_{i\,\psi(\cdot)}\rangle\,\widetilde{X_i}$. 
Due to (1), (2), we have (3). 
\end{pf}
\vskip0.6cm\par
\begin{lem}
Under the same assumption of Lemma 2.1, we have, 
\begin{equation}
\overline{\overline{\nabla}}_X(\overline{\overline{\nabla}}_YV)
=\sum_{i=1}^n
X_x(Y
\langle W,X_{i\,\psi(\cdot)}
\rangle)\,\pi_{\ast}(X_{i\,\psi(x)})
\in T_{\varphi(x)}(G/K), 
\end{equation}
at each $x\in M$, 
for every 
$C^{\infty}$ vector fields $X$ and $Y$ on $M$. 
\end{lem}
\begin{pf}
Let 
$Z:=\overline{\overline{\nabla}}_YV\in \Gamma(\varphi^{-1}T(G/K))$. 
Then, by Lemma 2.1 (1), we have 
\begin{equation}
\overline{\overline{\nabla}}_X(\overline{\overline{\nabla}}_YV)
=\overline{\overline{\nabla}}_XZ
=\pi_{\ast}\nabla^k_{(\psi_{\ast}X)^H}Z^H
\end{equation}
where by Lemma 2.1 (3), we have for every $y\in M$, 
$$
Z^H(y) 
=\sum_{i=1}^nY_y\langle W,X_{i\,\psi(\cdot)}\rangle\, 
X_{i\,\psi(y)},
$$
$Z(y)=\pi_{\ast}Z^H(y)\in T_{\varphi(y)}(G/K)$ and 
$Z\in \Gamma(\varphi^{-1}T(G/K))$. 
Then, at each $x\in M$, 
the right hand side of (2.17) which belong to $T_{\varphi(x)}(G/K)$,
coincides with the following: 
\begin{align*}
&
\sum_{j=1}^nX_x
\langle Z^H,X_{j\,\psi(\cdot)}\rangle\,\pi_{\ast}(X_{j\,\psi(x)})\\
&=\sum_{j=1}^nX_x
\,\langle
\sum_{i=1}^nY_{\bullet}
\langle W,X_{i\,\psi}\rangle\,X_{i\,\psi(\cdot)},
X_{j\,\psi(\cdot)}\rangle\,
\pi_{\ast}(X_{j\,\psi(x)})\\
&=\sum_{i,j=1}^n
X_x\,
\left(Y_{\bullet}\,\langle W,X_{i\,\psi}\rangle\right)
\,\delta_{ij}\,\pi_{\ast}(X_{j\,\psi(x)})\\
&=
\sum_{i=1}^nX_x\,\left(Y_{\bullet}\,\langle W,X_{i\,\psi}\rangle\right)
\,
\pi_{\ast}(X_{j\,\psi(x)}).
\end{align*}
Thus, we have (2.16). 
\end{pf}
\vskip0.6cm\par
\begin{prop}
The rough Laplacian 
$\overline{\Delta}$ acting on $\Gamma(\varphi^{-1}T(G/K))$ 
can be calculated as follows: 
For $V\in \Gamma(\varphi^{-1}T(G/K))$ with  
$V=\pi_{\ast}W$ for  
$W\in \Gamma(\psi^{-1}TG)$, 
\begin{equation}
(\overline{\Delta}V)(x)
=\sum_{i=1}^n
\Delta_x\,\langle W,X_{i\,\psi(\cdot)}\rangle \,
\pi_{\ast}(X_{i\,\psi(x)})
\in T_{\varphi(x)}(G/K), 
\end{equation}
for each $x\in M$. Here, since 
$f:\,M\ni x\mapsto \langle W(x),X_{i\,\psi(x)}\rangle_{\psi(x)}\in {\mathbb R}$ 
is a (local) $C^{\infty}$ function on $M$, the Laplacian 
$\Delta_x=\delta\,d$  acting on 
$C^{\infty}(M)$ works on $f$. 
\end{prop}
Indeed, if we recall the definition (2.3) of the rough Laplacian 
$\overline{\Delta}$, and due to Lemmas 2.1 and 2.2, we have 
\begin{align*}
\overline{\Delta}V
&=-\sum_{j=1}^m
\{
\overline{\overline{\nabla}}_{e_j}(\overline{\overline{\nabla}}_{e_j}V)
-\overline{\overline{\nabla}}_{\nabla_{e_j}e_j}V
\},\\
&=
-\sum_{i=1}^n
\sum_{j=1}^n(e_j{}^2-\nabla_{e_j}e_j)
\,\langle W,X_{i\,\psi(\cdot)}\rangle \,
\pi_{\ast}(X_{i\,\psi(x)})\\
&=
\sum_{i=1}^n
\Delta_x\,\langle W,X_{i\,\psi(\cdot)}\rangle \,
\pi_{\ast}(X_{i\,\psi(x)}).
\end{align*}
We have Proposition 2.3.\qed
\vskip0.6cm\par
\section{Determination of the bitension field}
Now, let $\theta$ be the Maurer-Cartan form on $G$, i.e., 
a $\frak g$-valued left invariant $1$-form on $G$ 
which is 
defined by 
$\theta_y(Z_y)=Z$
($y\in G$, 
$Z\in \frak g$). 
For every 
$C^{\infty}$ map $\varphi$ of $(M,g)$ into $(G/K,h)$ 
with a lift $\psi:\,M\rightarrow G$, 
let us consider 
a $\frak g$-valued $1$-form $\alpha$ on $M$ 
given by 
$\alpha=\psi^{\ast}\theta$ and 
the decomposition 
\begin{equation}
\alpha=\alpha_{\frak k}+\alpha_{\frak m}
\end{equation} 
corresponding to the decomposition 
${\frak g}={\frak k}\oplus{\frak m}$. 
Then, it is well known (see for example, 
\cite{DSU}) that 
\begin{lem} 
For every $C^{\infty}$ map 
$\varphi:\,(M,g)\rightarrow (G/K,h)$, 
\begin{equation}
t_{\psi(x)^{-1}\ast}
\tau(\varphi)=-\delta(\alpha_{\frak m}j+\sum_{i=1}^m
[\alpha_{\frak k}(e_i),\alpha_{\frak m}(e_i)],\quad (x\in M),
\end{equation}
where 
$\alpha=\varphi^{\ast}\theta$, and $\theta$ is the Maurer-Cartan form of $G$, 
$\delta(\alpha_{\frak m})$ is the co-differentiation of ${\frak m}$-valued 
$1$-form ${\alpha}_{\frak m}$ on $(M,g)$. 
\par
Thus, $\varphi:(M,g)\rightarrow (G/K,h)$ is harmonic if and only if 
\begin{equation}
-\delta(\alpha_{\frak m})+\sum_{i=1}^m
[\alpha_{\frak k}(e_i),\alpha_{\frak m}(e_i)]=0. 
\end{equation}
\end{lem}
\vskip0.6cm\par
Furthermore, we obtain 
\begin{th}
We have 
\begin{align}
&t_{\psi(x)^{-1}\ast}\tau_2(\varphi)
=\Delta_g
\left(
-\delta(\alpha_{\frak m})+\sum_{i=1}^m
\left[
\alpha_{\frak k}(e_i),\alpha_{\frak m}(e_i)
\right]
\right)
\nonumber\\
&
+\sum_{s=1}^m
\left[\left[
-\delta(\alpha_{\frak m})+\sum_{i=1}^m
\left[
\alpha_{\frak k}(e_i),\alpha_{\frak m}(e_i)
\right],\alpha_{\frak m}(e_s)
\right], 
\alpha_{\frak m}(e_s)
\right],
\end{align}
where 
$\Delta_g$ is the (positive) Laplacian of $(M,g)$ acting on $C^{\infty}$ 
functions on $M$, and 
$\{e_i\}_{i=1}^m$ is a local orthonormal frame field on $(M,g)$. 
\end{th}
\vskip0.6cm\par
Therefore, we obtain immediately the following two corollaries. 
\begin{cor}
Let $(G/K,h)$ be a Riemannian symmetric space, 
and $\varphi:\,(M,g)\rightarrow (G/K,h)$, a $C^{\infty}$ mapping. 
Then, we have: 
\par
$(1)$ \,\, the map $\varphi:\,(M,g)\rightarrow (G/K,h)$ is 
harmonic if and only if 
\begin{equation}
-\delta(\alpha_{\frak m})+\sum_{i=1}^m
\left[
\alpha_{\frak k}(e_i),\alpha_{\frak m}(e_i)
\right]
=0.
\end{equation}
\par $(2)$ \,\, The map $\varphi:\,(M,g)\rightarrow (G/K,h)$ is 
biharmonic if and only if 
\begin{align}
&
\Delta_g
\left(
-\delta(\alpha_{\frak m})+\sum_{i=1}^m
\left[
\alpha_{\frak k}(e_i),\alpha_{\frak m}(e_i)
\right]
\right)
\nonumber\\
&
+\sum_{s=1}^m
\left[\left[
-\delta(\alpha_{\frak m})+\sum_{i=1}^m
\left[
\alpha_{\frak k}(e_i),\alpha_{\frak m}(e_i)
\right],\alpha_{\frak m}(e_s)
\right], 
\alpha_{\frak m}(e_s)
\right]=0.
\end{align}
\end{cor}
\vskip0.6cm\par
\begin{cor}
Let $(G/K,h)$ be a Riemannian symmetric space, and 
$\varphi:\,(M,g)\rightarrow (G/K,h)$, a $C^{\infty}$ mapping with 
a {\em horizontal} lift 
$\psi:\,M\rightarrow G$, i.e., 
$\varphi=\pi\,\circ\,\psi$ and 
$\psi_x(T_xM)\subset H_{\psi(x)}$ which is equivalent to 
$\alpha_{\frak k}\equiv 0$. 
\par Then, we have: 
\par
$(1)$ the map $\varphi:\,(M,g)\rightarrow (G/K,h)$ is harmonic if and only if 
\begin{equation}
\delta(\alpha_{\frak m})=0,
\end{equation}
\par $(2)$ and 
the map $\varphi:\,(M,g)\rightarrow (G/K,h)$ is biharmonic if and only if 
\begin{equation}
\delta\,d\,\delta(\alpha_{\frak m})
+\sum_{s=1}^m
\left[\left[
\delta(\alpha_{\frak m}),\alpha_{\frak m}(e_s)
\right],\alpha_{\frak m}(e_s)
\right]=0.
\end{equation}
\end{cor}
\vskip0.6cm\par
{\it Proof of Theorem 3.2.} 
\par 
We need the following lemma: 
\begin{lem} 
The tension field $\tau(\varphi)$ of 
a $C^{\infty}$ map 
$\varphi:\,(M,g)\rightarrow (G/K,h)$ 
can be expressed as
$$\tau(\varphi)
=\pi_{\ast}W=\pi_{\ast}(W^H),$$
where $W\in \Gamma(\psi^{-1}TG)$,  and 
$W^H$ is the horizontal component of $W$ in the decomposition 
$W(x)=W^V(x)+W^H(x)\in T_{\psi(x)}G=V_{\psi(x)}\oplus H_{\psi(x)}$ 
$(x\in M)$. 
If we define an $\frak m$-valued function $\beta$ on $M$ by 
\begin{equation}
\beta:=\sum_{i=1}^n\langle W,\widetilde{X_i}\rangle X_i
=\sum_{i=1}^n\langle W^H,\widetilde{X_i}\rangle X_i,
\end{equation}
then, we have 
\begin{equation}
t_{\psi(x)\,\ast}{}^{-1}\tau(\varphi)=\pi_{\ast}\beta.
\end{equation}
\par
If we define $n$ $\frak m$-valued functions $\beta_i$ $(i=1,\cdots,n)$ on $M$ by
\begin{equation}
\beta_i:=\sum_{j=1}^n
\langle \psi_{\ast}e_i,X_{j\,\psi(\cdot)}\rangle\,X_j\in {\frak m}.
\end{equation}
Then, it holds that 
\begin{equation}
t_{\psi(x)\,\ast}{}^{-1}\varphi_{\ast}e_i=\pi_{\ast}\beta_i\,\, \text{and} \,\,
\beta_i=\alpha_{\frak m}(e_i),
\end{equation}
where $\alpha_{\frak m}$ is the $\frak m$-component of 
$\alpha:=\psi^{\ast}\theta$, and the Maurer-Cartan form on $G$.  
\end{lem}
Indeed, (3.10) and the first part of (3.12) follow from the definition of 
$\beta$ and 
the fact that 
\begin{align*}
\alpha(e_i)&=(\psi^{\ast}\theta)(e_i)\\
&=\theta(\psi_{\ast}e_i)\\
&=\theta\left(
\sum_{j=1}^n\langle\psi_{\ast}e_i,X_{j\,\psi(x)}\rangle\,X_{j\,\psi(x)}
+\sum_{j=n+1}^{\ell}\langle\psi_{\ast}e_i,X_{j\,\psi(x)}\rangle\,X_{j\,\psi(x)}
\right)\\
&=
\sum_{j=1}^n\langle\psi_{\ast}e_i,X_{j\,\psi(x)}\rangle\,
X_j
+\sum_{j=n+1}^{\ell}\langle\psi_{\ast}e_i,X_{j\,\psi(x)}\rangle\,
X_j\\
&\in {\frak m}\oplus {\frak k},
\end{align*} 
since $\alpha=\psi^{\ast}\theta$. Thus, we have 
$\beta_i=\alpha_{\frak m}(e_i)$. 
\qed
\vskip0.6cm\par
({\em Continued the proof of Theorem 3.2})
\quad
We have 
\begin{equation}
t_{\psi(x)\,\ast}{}^{-1}\varphi_{\ast}e_i
=\sum_{j=1}^n
\langle \psi_{\ast}e_i,X_{j\,\psi(x)}\rangle\,\pi_{\ast}(X_j)\in 
T_{o}(G/K),
\end{equation}
where $o=\{K\}\in G/K$ is the origin of $G/K$. 
Because, 
\begin{align*}
t_{\psi(x)\,\ast}{}^{-1}\varphi_{\ast}e_i=
t_{\psi(x)\,\ast}{}^{-1}\pi_{\ast}\,\psi_{\ast}e_i
=\pi_{\ast}\,L_{\psi(x)\,\ast}{}^{-1}\,\psi_{\ast}e_i
&=\pi_{\ast}(L_{\psi(x)\,\ast}{}^{-1}\psi_{\ast}\,e_i)_{\frak m}
\end{align*}
which coincides with 
\begin{equation*}
\sum_{j=1}^n\langle 
L_{\psi(x)\,\ast}{}^{-1}\,\psi_{\ast}e_i,X_j\rangle\,\pi_{\ast}(X_j)
=\sum_{j=1}^n
\langle \psi_{\ast}e_i,X_{j\,\psi(x)}\rangle\,\pi_{\ast}(X_j),
\end{equation*}
which imply (3.13). 
\par Thus, we have
\begin{equation}
\varphi_{\ast}e_i=\pi_{\ast}W_i\quad (i=1,\cdots,m),
\end{equation}
where $W_i\in \Gamma(\psi{}^{-1}TG)$ and 
${\frak m}$-valued functions 
$\widetilde{W}_i$ on $M$ 
$(i=1,\cdots,m)$ 
are given by 
\begin{align}
W_i(x)&:=\sum_{j=1}^n\langle \psi_{\ast}e_i,X_{j\,\psi(x)}\rangle\,X_{j\,\psi(x)},\\
\widetilde{W_i}(x)&:=
\sum_{j=1}^n\langle \psi_{\ast}e_i,X_{j\,\psi(x)}\rangle\,X_{j}\in {\frak m},
\end{align}
for each $x\in M$. 
\par
On the other hand, we have 
\begin{equation}
\tau(\varphi)=\pi_{\ast}W,
\end{equation}
where $W\in \Gamma(\psi^{-1}TG)$ and 
an 
$\frak m$-valued function $\widetilde{W}$ on $M$ are given by 
\begin{align}
W(x)&:=t_{\psi(x)\,\ast}\left(
-\delta(\alpha_{\frak m})+\sum_{i=1}^m
[\alpha_{\frak k}(e_i),\alpha_{\frak m}(e_i)]
\right),\\
\widetilde{W}&:=
-\delta(\alpha_{\frak m})+\sum_{i=1}^m[\alpha_{\frak k}(e_i),\alpha_{\frak m}(e_i)]
\end{align}
for each $x\in M$. 
And we also have 
\begin{equation}
t_{\psi(x)\,\ast}{}^{-1}\,\overline{\Delta}\tau(\varphi)(x)
=\Delta
\left(
-\delta(\alpha_{\frak m})+\sum_{i=1}^m
[\alpha_{\frak k}(e_i),\alpha_{\frak m}(e_i)]
\right)(x)\quad (x\in M),
\end{equation}
where 
$\Delta=\delta\,d$ is the positive Laplacian acting on the space of all 
$C^{\infty}$ ${\frak m}$-valued functions on $M$. 
\par
We want to calculate 
${\mathcal R}(\tau(\varphi))=\sum_{i=1}^m
R^h(\tau(\varphi),\varphi_{\ast}e_i)\varphi_{\ast}e_i$. 
Indeed, we have 
\begin{align}
t_{\psi(x)\,\ast}{}^{-1}\,
{\mathcal R}(\tau(\varphi))
&=-\sum_{i=1}^m
[[
\widetilde{W},\widetilde{W_i}],\widetilde{W_i}]\nonumber\\
&=
-\sum_{s=1}^m
[[
-\delta(\alpha_{\frak m})+\sum_{i=1}^m
[\alpha_{\frak k}(e_i),\alpha_{\frak m}(e_i)],\nonumber\\
&\qquad \qquad\qquad\quad
\alpha_{\frak m}(e_s)],
\alpha_{\frak m}(e_s)].
\end{align}
Here, we used the formula of the curvature $R^h$ of the Riemannian symmetric 
space $(G/K,h)$ (\cite{KN}, p. 202, p.231, Theorem 3.2) : 
\begin{equation*}
(R^h(X,Y)Z)_o
=-[[X,Y],Z]_o
\quad (X,Y,Z\in {\frak m}). 
\end{equation*}
Thus, we obtain Theorem 3.2. \qed
\vskip0.6cm\par
Let us recall the {\em integrability condition} 
for a $C^{\infty}$ mapping $\varphi:\,(M,g)\rightarrow 
(G/K,h)$.  
The Maurer-Cartan form 
$\theta$ on $G$ satisfies 
\begin{equation}
d\theta+\frac12[\theta\wedge\theta]=0, 
\end{equation}
so that the pull back $\alpha=\psi^{\ast}\theta$ 
of $\theta$ by the lift 
$\psi:\,M\rightarrow G$ of 
$\varphi:\,M\rightarrow G/K$ also satisfies that 
\begin{equation}
d\alpha+\frac12[\alpha\wedge\alpha]=0,
\end{equation} 
which is equivalent to 
\begin{equation}
\left\{
\begin{aligned}
&d\alpha_{\frak k}+\frac12[\alpha_{\frak k}\wedge\alpha_{\frak k}]
+\frac12[\alpha_{\frak m}\wedge\alpha_{\frak m}]=0,\\
&d\alpha_{\frak m}+[\alpha_{\frak k}\wedge\alpha_{\frak m}]=0.
\end{aligned}
\right.
\end{equation}
\par
Summarizing up the above, we have 
\begin{th}
Let $(M,g)$ be an $m$-dimensional compact Riemannian manifold, 
$(G/K,h)$, an $n$-dimensional Riemannian symmetric space, 
$\pi:\,G\rightarrow G/K$, the projection, and 
$\varphi:\,(M,g)\rightarrow (G/K,h)$, a $C^{\infty}$ mapping with 
a local lift 
$\psi:\,M\rightarrow G$, 
$\varphi=\pi\circ\psi$. Let $\alpha=\psi^{\ast}\theta$ be the pull back of the Maurer-Cartan form $\theta$, and $\alpha=\alpha_{\frak k}+\alpha_{\frak m}$, the decomposition of $\alpha$ 
corresponding to the Cartan decomposition 
${\frak g}={\frak k}\oplus{\frak m}$. 
\par
$(I)$ The mapping $\varphi:\,(M,g)\rightarrow (G/K,h)$ is {\em harmonic} 
if and only if 
\begin{equation}
-\delta(\alpha_{\frak m})+\sum_{i=1}^m[\alpha_{\frak k}(e_i),\alpha_{\frak m}(e_i)]=0, 
\end{equation}
where $\delta$ is the co-differentiation, and $\{e_i\}_{i=1}^m$ is a local orthonormal frame field on $(M,g)$. 
\par
Furthermore, $\varphi:\,(M,g)\rightarrow (G/K,h)$ is {\em biharmonic} 
if and only if 
\begin{align}
&\Delta\left(-\delta(\alpha_{\frak m})+\sum_{i=1}^m[\alpha_{\frak k}(e_i),\alpha_{\frak m}(e_i)]\right)\nonumber\\
&\,\,+\sum_{s=1}^m
\left[\left[
-\delta(\alpha_{\frak m})+\sum_{i=1}^m[\alpha_{\frak k}(e_i),\alpha_{\frak m}(e_i)],\alpha_{\frak m}(e_s)
\right],
\alpha_{\frak m}(e_s)\right]
=0, 
\end{align}
where $\Delta=\delta d$ is the (positive) Laplacian of $(M,g)$ 
acting on the space of $\frak g$-valued $C^{\infty}$ functions on $(M,g)$. 
\par
$(II)$ Conversely, let $\alpha=\alpha_{\frak k}+\alpha_{\frak m}$ 
be a $\frak g$-valued $1$-form on $(M,g)$. 
If $\alpha$ satisfies $(3.23)$ or $(3.24)$, and satisfies $(3.25)$ 
(resp. $(3.26)$), 
then, there exists a $C^{\infty}$-mapping $\varphi$ of $M$ into $G$ with 
a local lift $\psi:\,M\rightarrow G$, $\varphi=\pi\circ\psi$ and 
the initial value $\varphi(p)=a\in G$ at some $p\in M$ 
such that 
$\alpha=\psi^{\ast}\theta$ and 
$\varphi$ is a harmonic (resp. biharmonic) map of $(M,g)$ into $(G/K,h)$. 
\end{th}
\vskip0.6cm\par
\section{Biharmonic curves into Riemannian symmetric spaces}
\subsection
\quad 
Let $\varphi:\,({\mathbb R},g_0)\rightarrow (G/K,h)$ 
be a $C^{\infty}$ curve, 
and $\psi:\,{\mathbb R}\rightarrow G$, a lift of $\varphi$, 
$(\varphi=\pi\circ\psi)$.  
Then, 
$\alpha=\psi^{\ast}\theta=\psi^{-1}d\psi=F(t)dt$ is a 
${\frak g}$-valued $1$-form on $\mathbb R$ and 
$F$ is a $\frak g$-valued function on $\mathbb R$ satisfying 
$\psi(t)^{-1}\frac{d\psi}{dt}=F(t)$.
Conversely, for a $\frak g$-valued $C^{\infty}$ function $F(t)$ on 
$\mathbb R$,
there exists a 
unique $C^{\infty}$-curve 
$\psi:\,{\mathbb R}\rightarrow G$ which satisfies that 
\begin{equation}
\left\{
\begin{aligned}
\psi(t)^{-1}\frac{d\psi}{dt}=F(t),\\
\psi(0)=x\in G. 
\end{aligned}
\right.
\end{equation}
\par
To give an explicit solution $\psi$ of  $(4.1)$ is very difficult for us, in general, since $G$ is not abelian. 
However, corresponding to the decomposition 
${\frak g}={\frak k}\oplus{\frak m}$, we decompose 
$F(t)=F_{\frak k}(t)+F_{\frak m}(t)$, 
$\alpha_{\frak k}=F_{\frak k}(t)dt$, 
and $\alpha_{\frak m}=F_{\frak m}(t)dt$, 
so we have
$$\delta\alpha=-(\overline{\nabla}_{e_1})(\alpha(e_1))=-\nabla^h_{e_1}(\alpha(e_1))=-e_1(F(t))=-F'(t),$$
and 
$$\delta\alpha_{\frak m}=-F_{\frak m}{}'(t).
$$
Thus the harmonic map equation $(3.25)$ is 
\begin{equation}
F_{\frak m}{}'(t)+[F_{\frak k}(t),F_{\frak m}(t)]=0,
\end{equation}
and the biharmonic map equation $(3.26)$ is 
\begin{align}
&-\frac{d^2}{dt^2}\left(
F_{\frak m}{}'(t)+[F_{\frak k}(t),F_{\frak m}(t)]
\right)\nonumber\\
&\qquad+\left[\left[
F_{\frak m}{}'(t)+[F_{\frak k}(t),F_{\frak m}(t)],F_{\frak m}
\right],F_{\frak m}\right]=0.
\end{align}
In these cases, 
the integrability condition $(3.23)$ always holds, so that the existence of $\psi$ of $(4.1)$ is always true. 
\par
Let us recall that a lift $\psi(t)$ is {\em horizontal} if 
$\psi_{\ast}(T_xM)\subset L_{\ast\,\psi(x)}({\frak m})$ 
if and only if 
$F_{\frak k}\equiv0$. 
In this case, 
$(4.2)$ is equivalent to 
\begin{equation}
F_{\frak m}{}'(t)=0, 
\end{equation}
which implies that 
$F_{\frak m}(t)=X\in {\frak m}$ (constant). So that 
$F(t)=X\in {\frak m}$. 
Then, we have 
\begin{equation}
\psi(t)=x\,\exp(tX), \quad \varphi(t)=x\exp(tX)\,K\in G/K. 
\end{equation}
Furthermore,  $(4.3)$ is equivalent to 
\begin{equation}
-F_{\frak m}{}'''(t)+[[F_{\frak m}{}'(t),F_{\frak m}(t)],F_{\frak m}(t)]=0. 
\end{equation}
\vskip0.6cm\par
{\it Example 4.1.} \quad Assume that 
$(G/K,h)$ is {\em of the Euclidean type}. In this case, $\frak m$ is an abelian ideal and $\frak k$ acts on $\frak m$ by 
$[T,X]=T\cdot X$ $(T\in {\frak k}, X\in {\frak m})$ regarding 
$\frak k$ as a subalgebra of ${\frak gl}({\frak m})$. Then, we have 
\par
$(1)$ $\varphi:\,({\mathbb R},g_0)\rightarrow (G/K,h)$ is harmonic 
if and only if 
\begin{equation}F_{\frak m}{}'(t)+F_{\frak k}(t)\cdot F_{\frak m}(t)=0.
\end{equation}
\par
$(2)$ $\varphi:\,({\mathbb R},g_0)\rightarrow (G/K,h)$ is biharmonic 
if and only if 
\begin{equation}
\frac{d^2}{dt^2}\left(F_{\frak m}{}'(t)+F_{\frak k}(t)\cdot F_{\frak m}(t)\right)=0
\end{equation}
which is equivalent to 
\begin{equation}
F_{\frak m}{}'(t)+F_{\frak k}(t)\cdot F_{\frak m}(t)=At+B
\end{equation}
for some $A$ and $B$ in ${\frak m}$. 
Thus, if $\psi:\,({\mathbb R},g_0)\rightarrow G$ is horizontal, i.e., 
$F_{\frak k}\equiv0$, then, 
$F_{\frak m}(t)=C$ (a constant vector in ${\frak m}$) 
for the case $(1)$, and 
$F_{\frak m}(t)=At^2+Bt+C$ for the case $(2)$. 
If $[A,B]=[B,C]=[C,A]=0$, 
then 
$\psi(t)=\exp(t^2\,A+t\,B+C)$ and 
$\varphi(t)=\psi(t)\cdot\{K\}$ is a biharmonic curve in a Riemannian symmetric space $(G/K,h)$ of the Euclidean type.  
\vskip0.6cm\par
\subsection{Biharmonic curves into rank one symmetric spaces} \quad 
In this subsection, we study biharmonic curves in a compact symmetric spaces $(G/K,h)$. 
\vskip0.3cm\par
{\it $(1)$ Case of the unit sphere $(S^n,h)$.} \quad 
Let $G=SO(n+1)$ act on ${\mathbb R}^{n+1}$ linearly, 
and $K=SO(n)$ be the isotropy subgroup of $G$ at the origin 
$o={}^{\rm t}(1,0,\cdots,0)$. Their Lie algebras 
${\frak g}={\frak so}(n+1)$, ${\frak k}={\frak so}(n)$ and the Cartan 
decomposition 
${\frak g}={\frak k}\oplus {\frak m}$ are given by 
\begin{align*}
{\frak g}&={\frak so}(n+1)=\{X\in {\frak gl}(n+1):\,X+{}^{\rm t\!}X=O\},\\
{\frak k}&={\frak so}(n)=
\left\{
\begin{pmatrix}
0&\strut\vrule&0&\cdots&0\\
\noalign{\vskip-1pt}
\multispan6\hrulefill\\
0&\strut\vrule& & &\\
\vdots&\strut\vrule& &X_1&\\
0&\strut\vrule& & & &
\end{pmatrix}
:\,X_1\in {\frak gl}(n),\,X_1+{}^{\rm t\!}X_1=O
\right\},
\\
{\frak m}&=\left\{
\begin{pmatrix}
0&\strut\vrule&-{}^{\rm t\!}u\\
\noalign{\vskip-1pt}
\multispan3\hrulefill\\
u&\strut\vrule& O\\
\end{pmatrix}
:\,u={}^{\rm t\!}(u_1,\cdots,u_n)\in {\mathbb R}^n
\right\}.
\end{align*}
For a ${\frak m}$-valued $C^{\infty}$  function $F_{\frak m}(t)$ given by 
\begin{equation}
F_{\frak m}(t)=
\begin{pmatrix}
0&\strut\vrule&-u_1(t)&\cdots&-u_n(t)\\
\noalign{\vskip-1pt}
\multispan6\hrulefill\\
u_1(t)&\strut\vrule& & &\\
\vdots&\strut\vrule& & O &\\
u_n(t)&\strut\vrule& & & &
\end{pmatrix},
\end{equation} 
and $F_{\frak k}\equiv0$, 
the biharmonic map equation $(4.7)$ 
is equivalent to 
\begin{equation}
-u_i{}'''+\sum_{j=1}^n(u_i\,u_j{}'-u_i{}'\,u_j)u_j=0
\quad (i=1,\cdots,n)
\end{equation}
which is also equivalent to 
\begin{equation}
-u'''+\langle u',u\rangle u-\langle u,u\rangle u'=0,
\end{equation}
where the inner product $\langle\,,\,\rangle$ 
on ${\mathbb R}^n$ is given by 
$\langle u,v\rangle=\sum_{i=1}^nu_iv_i$ for $u,v\in {\mathbb R}^n$. 
\vskip0.3cm\par
\underline{{\it Case of $n=2$}}.  
Our problem is to find a $C^{\infty}$ 
plane curve which satisfies $(4.12)$. 
To do it, we assume that 
$u(t)$ is reparametrized in such a way that 
$u(s)$ is a tangent curve of a plane curve 
${\mathbf p}(s)$: 
$u(s)={\mathbf p}'(s)={\mathbf e}_1(s)$. 
For the other cases, we have no idea to solve $(4.12)$.   
Recall the Frenet-Serret formula for a plane curve ${\mathbf p}(s)$:
\begin{equation}
\left\{
\begin{aligned}
{\mathbf p}'(s)&={\mathbf e}_1(s),\\
{\mathbf e}_1'(s)&=\kappa(s)\,{\mathbf e}_2(s),\\
{\mathbf e}_2'(s)&=-\kappa(s)\,{\mathbf e}_1(s). 
\end{aligned}
\right.
\end{equation}
Now we have 
\begin{align}
u&={\mathbf e}_1,\\
u'&={\mathbf e}_1'=\kappa\,{\mathbf e}_2,\\
u''&=\kappa'\,{\mathbf e}_2+\kappa\,{\mathbf e}_2'=
-\kappa^2\,{\mathbf e}_1+\kappa'\,{\mathbf e}_2,\\
u'''&=-3\kappa\,\kappa'\,{\mathbf e}_1+(\kappa''-\kappa^3)\,{\mathbf e}_2.
\end{align}
Since $\langle u',u\rangle=0$ and $\langle u,u\rangle=1$, 
$(4.12)$ is equivalent to 
\begin{align}
-3\kappa \,\kappa'&=0,\\
\kappa''-\kappa^3&=-\kappa,
\end{align}
By $(4.18)$, $\kappa=c$ (a constant), and by $(4.19)$, 
$c=0,\,1, \,-1$. Thus, we have \par
$(i)$ In the case of $c=0$, 
\begin{equation}
{\mathbf p}(s)=s\,{\mathbf a}+{\mathbf b}, \,\,
u(s)={\mathbf a},\,\,({\mathbf a, b}\in {\mathbf R}^2), 
\end{equation}
\par
$(ii)$ in the case of $c=1$,
\begin{equation}
{\mathbf p}(s)=(\cos s,\sin s),\,\,
u(s)=(-\sin s,\cos s),
\end{equation}
\par
$(iii)$ in the case of $c=-1$, 
\begin{equation}
{\mathbf p}(s)=(\cos s,-\sin s),\,\,
u(s)=(-\sin s,-\cos s). 
\end{equation}
Now it is easy to find $\psi:{\mathbb R}\rightarrow G$ and 
$\varphi(t)=\psi(t)\,\{K\}\in G/K$ satisfying 
$\psi(t)^{-1}\frac{d\psi}{dt}=F(t)=F_{\frak m}(t)$ for such $u(t)$ in $(4.12)$. 
\par
Case $(i)$: \quad 
If ${\mathbf a}={}^{\rm t\!}(a,b)\in {\mathbb R}^2$,  we have due to $(4.1)$, 
\begin{equation}
\varphi(t)=\psi(t)\,\{K\}=x\begin{pmatrix}
\cos (t\,\sqrt{a^2+b^2})\\
\frac{a}{\sqrt{a^2+b^2}}\sin (t\,\sqrt{a^2+b^2})\\
\frac{b}{\sqrt{a^2+b^2}}\sin (t\,\sqrt{a^2+b^2})
\end{pmatrix},
\end{equation}
which is a great circle of the standard $2$-sphere $(S^2,h)$. 
\par
Cases $(ii)$ and $(iii)$:  In these cases, 
if we assume $F_{\frak k}\equiv0$, 
we have 
\begin{equation}
F_{\frak m}(t)=
\begin{pmatrix}
0&\strut\vrule&\sin t& -\cos t\\
\noalign{\vskip-1pt}
\multispan5\hrulefill\\
-\sin t&\strut\vrule&0&0 &\\
\cos t&\strut\vrule& 0& 0& 
\end{pmatrix},
\end{equation}
for Case $(ii)$, and 
\begin{equation}
F_{\frak m}(t)=
\begin{pmatrix}
0&\strut\vrule&\sin t& \cos t\\
\noalign{\vskip-1pt}
\multispan5\hrulefill\\
-\sin t&\strut\vrule&0&0 &\\
-\cos t&\strut\vrule& 0& 0& 
\end{pmatrix},
\end{equation}
for Case $(iii)$. 
In these cases, because of 
$[F_{\frak m}(t),F_{\frak m}{}'(t)]\not=0$, 
it is difficult for us to give explicitly a unique solution 
of the initial value problem of 
\begin{equation}\psi(t)^{-1}\frac{\psi(t)}{dt}=F(t) \quad\text{and}\,\, \psi(0)=a\in SO(3).
\end{equation}
\vskip0.3cm\par
\underline{\it Case of $n=3$}. 
In this case, we have to solve 
for a $C^{\infty}$ curve 
$u:{\mathbb R}\rightarrow {\mathbb R}^3$, 
the equation $(4.12)$ 
which is equivalent to
\begin{equation}
-u'''+u\times(u\times u')=0. 
\end{equation} 
To do it, 
we assume that $u(t)$ is parametrized in such a way 
that $u(s)$ is a tangent curve of a $C^{\infty}$ curve in 
${\mathbb R}^3$, 
${\mathbf p} (s):\,u(s)={\mathbf p}'(s)={\mathbf e}_1(s)$. 
Recall the Frene-Serret formula for a curve ${\mathbf p}(s)$: 
\begin{equation}
\left\{
\begin{aligned}
{\mathbf p}&={\mathbf e}_1\\
{\mathbf e}_1{}'&=\qquad \kappa{\mathbf e}_2\\
{\mathbf e}_2{}'&=-\kappa{\mathbf e}_1\quad +\tau{\mathbf e}_3\\
{\mathbf e}_3{}'&=\qquad -\tau{\mathbf e}_2
\end{aligned}
\right.
\end{equation}
where $\kappa$ and $\tau$ are the curvature and torsion 
of ${\mathbf p}(s)$, respectively. 
By making use of $(4.28)$, we have 
\begin{equation}
\left\{
\begin{aligned}
u'&=\kappa{\mathbf e}_2\\
u''&=-\kappa^2{\mathbf e}_1+\kappa'{\mathbf e}_2+\kappa\tau{\mathbf e}_3\\
u'''&=
-3\kappa\kappa'{\mathbf e}_1+
(\kappa''-\kappa^3-\kappa\tau^2){\mathbf e}_2
+(2\kappa'\tau+\kappa\tau'){\mathbf e}_3.
\end{aligned}
\right.
\end{equation}
Thus, $(4.29)$ is equivalent to 
\begin{equation}
\left\{
\begin{aligned}
-3\kappa\kappa'&=0\\
\kappa''-\kappa^3-\kappa\tau^2&=-\kappa\\
2\kappa'\tau+\kappa\tau'&=0.
\end{aligned}
\right.
\end{equation}
By the first equation of $(4.30)$, $\kappa=\kappa_0$ (a constant). 
In the case $\kappa_0=0$, $u(t)={\mathbf a}\in {\mathbb R}^3$ 
(a constant vector). 
In the case $\kappa_0\not=0$, by the third equation of $(4.30)$, 
$\tau=\tau_0$ (a constant). By the second equation of $(4.30)$, 
$\kappa_0{}^2+\tau_0{}^2=1$. Then, 
${\mathbf p}(s)={}^{\rm t\!}(a \cos t,a \sin t,b t)$, with 
$s=\sqrt{a^2+b^2}\,t$. 
Here, 
$\kappa_0=a/(a^2+b^2)$, and $\tau_0=b/(a^2+b^2)$, 
and 
$1=\kappa_0{}^2+\tau_0{}^2=1/(a^2+b^2)$, i.e., 
$a^2+b^2=1$. Therefore, we have 
\begin{equation}
\left\{
\begin{aligned}
{\mathbf p}(t)&={}^{\rm t\!}(a\cos t, a\sin t,bt), 
\\
u(t)&={\mathbf p}'(t)={}^{\rm t\!}(-a\sin t,a\cos t,b),
\end{aligned}
\right.
\end{equation}
where $a$ and $b$ are constants with $a^2+b^2=1$. 
Thus, $F_{\frak m}(t)$ with $F_{\frak k}\equiv 0$, is given by 
 \begin{equation}
F_{\frak m}(t)=
\begin{pmatrix}
0&\strut\vrule&a\sin t& -a\cos t&-b\\
\noalign{\vskip-1pt}
\multispan5\hrulefill\\
-a\sin t&\strut\vrule& & &\\
a\cos t&\strut\vrule& & O & \\
b&\strut\vrule& & &
\end{pmatrix}.
\end{equation}
Thus, we derived
only to solve explicitly the initial value problem $(4.26)$ 
 which is a difficult problem for us. 
\vskip0.3cm\par
\underline{Case of $n\geq 2$}.  
\quad 
In this case, the other-type solutions exist: \par
Let $u=(u_1,\cdots, u_n)=(0,\cdots,0,\overbrace{v}^{\text{$i$ th}},0,\cdots,0)$ $(i=1,\cdots,n)$.  Then, for such $u$, the equation 
$(4.12)$ is reduced to 
$v'''=0$. Thus, we have $v(t)=D_t:=at^2+bt+c$ for some constants 
$a$, $b$ and $c$. 
 Thus, $F_{\frak m}(t)$ 
 is given by 
 \begin{equation}
 F_{\frak m}(t)=D_t\,
 \begin{pmatrix}
 0&\strut\vrule&0&\cdots&-1&\cdots&0&\\
\noalign{\vskip-1pt}
\multispan8\hrulefill\\
0&\strut\vrule& & & & & &\\
\vdots&\strut\vrule& & & & & & \\
1&\strut\vrule& & & O & & &\\
\vdots&\strut\vrule& & & & & &\\
0&\strut\vrule & & & & & &
 \end{pmatrix}.
 \end{equation}
 Thus, 
 \begin{align}
 \psi(t)=x\,\exp
 \left(
 \int^t_0F(s)ds
 \right)
 =x\,
  \begin{pmatrix}
\cos d_t&\strut\vrule&0&\cdots&-\sin d_t&\cdots&0&\\
\noalign{\vskip-1pt}
\multispan8\hrulefill\\
0&\strut\vrule& 0& \cdots& 0&\cdots & 0&\\
\vdots&\strut\vrule& \vdots& & \vdots& &\vdots & \\
\sin d_t&\strut\vrule& 0&\cdots & \cos d_t & \cdots& 0&\\
\vdots&\strut\vrule&\vdots & & \vdots& & \vdots&\\
0&\strut\vrule &0 & \cdots& 0&\cdots & 0&
 \end{pmatrix},\nonumber
 \end{align}
 where $d_t:=\frac{a}{3}t^3+\frac{b}{2}t^2+ct$. 
 So, we have a {\em biharmonic curve} into $(S^n,h)$: 
 \begin{equation}
 \varphi(t)=\psi(t)\{K\}
 =x\,{}^{\rm t\!}(\cos d_t,0,\cdots,0,\sin d_t,0,\cdots,0),
 \end{equation}
 for $x\in SO(n+1)$, 
 where $d_t:=\frac{a}{3}t^3+\frac{b}{2}t^2+ct$. Furthermore, 
 $\varphi(t)$ is harmonic if and only if $a=b=0$. 
 \vskip0.3cm\par
 {\it $(2)$ Case of the complex projective space $({\mathbb C}P^n,h)$}. 
 \par 
 Let 
 $G=SU(n+1)$ act on the projective space linearly on 
 ${\mathbb C}P^n
 =\{[z]:\,z\in {\mathbb C}^{n+1}\backslash \{\mathbf 0\}\}$, and 
 $K$, the isotropy subgroup of $G$ 
 at $o={}^{\rm t\!}[1,0,\cdots,0]$. The Cartan decomposition 
 ${\frak g}={\frak k}\oplus {\frak m}$ is given by 
 \begin{align*}
 {\frak g}&=\{X\in {\frak gl}(n+1,{\mathbb C}):\,X+{}^{\rm t\!}\overline{X}=O,\,\,{\rm tr}X=0\},\\
 {\frak k}&=\left\{
 \begin{pmatrix}
 \sqrt{-1}\,a&0\\
 0&X
 \end{pmatrix}
 :\,a\in {\mathbb R}, X\in {\frak gl}(n,{\mathbb C}),\,
 {}^{\rm t\!}\overline{X}+X=O,\,\right.\\
 &
 \left.\qquad\qquad\qquad\qquad\qquad\qquad\qquad\qquad\sqrt{-1}a+{\rm tr}X=0
 \right\},\\
 {\frak m}&=\left\{
 \begin{pmatrix}
 0&-{}^{\rm t\!}\overline{z}\\
 z&O
 \end{pmatrix}:\,z\in {\mathbb C}^n
 \right\}. 
 \end{align*}
 For a $C^{\infty}$ ${\frak m}$-valued function $F_{\frak m}(t)$ 
 given by 
 \begin{equation}
 F_{\frak m}(t)
 =\begin{pmatrix}
0&\strut\vrule&-\,\overline{z_1(t)}&\cdots&-\,\overline{z_n(t)}\\
\noalign{\vskip-1pt}
\multispan6\hrulefill\\
z_1(t)&\strut\vrule& & &\\
\vdots&\strut\vrule& & O &\\
z_n(t)&\strut\vrule& & & &
\end{pmatrix},
 \end{equation}
 where $z_i(t)=u_i(t)+\sqrt{-1}v_i(t)$ 
 $u_i(t)$ and $v_i(t)$ are real valued $C^{\infty}$ functions 
 $(i=1,\cdots,n)$, and $F_{\frak k}\equiv 0$, 
 the biharmonic map equation $(4.6)$ is equivalent to 
 \begin{equation}
 -z_i{}'''+\sum_{j=1}^n
 \left\{
 (z_i\,\overline{z_j}'-z_i{}'\,\overline{z_j})\,z_j
 -z_i
(\overline{z_j}z_j{}'-\overline{z_j}'\,z_j)
 \right\}=0
 \end{equation}  
 for all $i=1,\cdots,n$. 
 Notice here that this 
 $(4.36)$ can be written as 
 \begin{equation}
 -z'''+2\,\langle z,z'\rangle \,z-\langle z',z\rangle\,z -\langle z,z\rangle\,z'=0,
 \end{equation}
 where $\langle z,w\rangle=\sum_{i=1}^nz_i\overline{w_i}$ 
 for two 
${\mathbb C}^n$-valued functions  $z$ and $w$ in $t$. 
 If we write $z=u+\sqrt{-1}v$, where 
 $u$ and $v$ are ${\mathbb R}^n$-valued functions, 
 then $(4.37)$ is equivalent to 
 \begin{equation}
 \left\{
 \begin{aligned}
 -u'''+4n\,(-v^2\,u'+u\,v\,v')&=0\\
 -v'''+4n\,(u\,v\,u'-u^2\,v')&=0.
 \end{aligned}
 \right.
 \end{equation}
 One can find the following solutions of $(4.38)$: 
 \vskip0.3cm\par
\qquad $(i)$ $u=D_t=a\,t^2+b\,t+c$ and $v\equiv0$,
 \par
\qquad $(ii)$ $u\equiv0$ and $v=D_t=a\,t^2+b\,t+c$, or 
 \par
\qquad $(iii)$ $u=v=D_t=a\,t^2+b\,t+c$, 
 \par\noindent
 where $a$, $b$ and $c$ are constant vectors in ${\mathbb R}^n$. 
 Corresponding to these, we can find $F_{\frak m}(t)$ of $(4.35)$
 as follows: 
\begin{equation}
 F_{\frak m}(t)
 =D_t
 \begin{pmatrix}
0&\strut\vrule&-\,\overline{z_1(t)}&\cdots&-\,\overline{z_n(t)}\\
\noalign{\vskip-1pt}
\multispan6\hrulefill\\
z_1(t)&\strut\vrule& & &\\
\vdots&\strut\vrule& & O &\\
z_n(t)&\strut\vrule& & & &
\end{pmatrix},
 \end{equation}
where $z_1(t),\cdots,z_n(t)$ are 
\par
Case $(i)$:\quad $z_1(t)=\cdots=z_n(t)=1$, 
\par
Case $(ii)$:\quad $z_1(t)=\cdots=z_n(t)=\sqrt{-1}$, 
\par
Case $(iii)$:\quad $z_1(t)=\cdots=z_n(t)=1+\sqrt{-1}$, 
\par\noindent
correspondingly. In each cases, we can find $\psi(t)$ by the same way as the case of $(S^n,h)$, and 
a {\it biharmonic curve} in 
$({\mathbb C}P^n,h)$: 
\par
Case $(i)$:\quad 
$\varphi(t)=x\,{}^{\rm t\!}[\cos(\sqrt{n}\,d_t),
\frac{1}{\sqrt{n}}\sin(\sqrt{n}\,d_t),\cdots,\frac{1}{\sqrt{n}}\sin(\sqrt{n}\,d_t)],$
\par
Case $(ii)$:\quad 
$\varphi(t)=x\,{}^{\rm t\!}[\cos(\sqrt{n}\,d_t),
\frac{\sqrt{-1}}{\sqrt{n}}\sin(\sqrt{n}\,d_t),\cdots,\frac{\sqrt{-1}}{\sqrt{n}}\sin(\sqrt{n}\,d_t)],$
\par
 Case $(iii)$:\par\qquad
$\varphi(t)=x\,{}^{\rm t\!}[\cos(\sqrt{2n}\,d_t),
\frac{1+\sqrt{-1}}{\sqrt{2n}}\sin(\sqrt{2n}\,d_t),\cdots,\frac{1+\sqrt{-1}}{\sqrt{2n}}\sin(\sqrt{2n}\,d_t)],$
where 
$d_t:=\frac{a}{3}\,t^3+\frac{b}{2}\,t^2+c\,t$, $a$, $b$ and $c$ 
are constant real numbers, and $x\in SU(n+1)$. 
Each 
$\varphi:\, ({\mathbb R},g_0)\rightarrow ({\mathbb C}P^n,h)$ 
is {\it harmonic} if and only if $a=b=0$. 
\vskip0.3cm\par
 \vskip0.3cm\par
 {\it $(3)$ Case of the quaternion projective space $({\mathbb H}P^n,h)$}. 
\par \quad Let 
 $G=Sp(n+1)=\{x\in U(2n+2)\vert\,
 {}^{\rm t\!}x\,J_{n+1}x=J_{n+1}\}$, 
 where $J_{n+1}=\begin{pmatrix} O&I_{n+1}\\-I_{n+1}&O \end{pmatrix}$, and 
 $I_{n+1}$ is the identity matrix of order $n+1$. 
 $G$ acts on the quaternion projective space linearly on 
 ${\mathbb H}P^n
 =\{[z]:\,z\in {\mathbb H}^{n+1}\backslash \{\mathbf 0\}\}$, and 
 $K=Sp(1)\times Sp(n)$ is the isotropy subgroup $K$ of $G$ 
 at $o={}^{\rm t\!}[1,0,\cdots,0]$. The Cartan decomposition 
 ${\frak g}={\frak k}\oplus {\frak m}$ is given by 
 \begin{align*}
 {\frak g}&={\frak sp}(n+1)=
 \left\{
 \begin{pmatrix}
 A&B\\
 -\overline{B}&\overline{A}
 \end{pmatrix}\vert\,
 A,B\in M_{n+1}({\mathbb C}),\,
 {}^{\rm t\!}\overline{A}+A=O, {}^{\rm t\!}B=B
 \right\},\\
 {\frak k}&={\frak sp}(1)\times {\frak sp}(n)=
 \left\{\begin{pmatrix}
x&0&\strut\vrule&y&0&\\
0&X&\strut\vrule&0&Y&\\
\noalign{\vskip-1pt}
\multispan6\hrulefill\\
-\overline{y}&0&\strut\vrule&\overline{x} & 0&\\
0&-\overline{Y}&\strut\vrule&0 & \overline{X} &
\end{pmatrix}
\vert\,x\in \sqrt{-1}{\mathbb R},\,y\in {\mathbb C}, \right.\\
&\left.\qquad\qquad\qquad\qquad\qquad\qquad\quad
X,Y\in M_n({\mathbb C}), 
{}^{\rm t\!}\overline{X}+X=0, {}^{\rm t\!}Y=Y
 \right\},\\
 {\frak m}&=
 \left\{
 \begin{pmatrix}
0&Z&\strut\vrule&0&W&\\
-{}^{\rm t\!}\overline{Z}&O&\strut\vrule&{}^{\rm t\!}W&O&\\
\noalign{\vskip-1pt}
\multispan6\hrulefill\\
\noalign{\vskip2pt}
0&-\overline{W}&\strut\vrule&0 & \overline{Z}&\\
-{}^{\rm t\!}\overline{W}&O&\strut\vrule&-{}^{\rm t\!}Z & O&
\end{pmatrix}
\vert\,
\,Z, W\in M(1,n,{\mathbb C})
 \right\}. 
 \end{align*}
 For a $C^{\infty}$ ${\frak m}$-valued function $F_{\frak m}(t)$ 
 given by 
 \begin{equation}
 F_{\frak m}(t)=
 \begin{pmatrix}
0
&Z
&0&W&\\
-{}^{\rm t\!}\overline{Z}
&O 
&{}^{\rm t\!}W&O&\\
0
&-\overline{W}
&0
&\overline{Z}&\\
-{}^{\rm t\!}\overline{W}& O
&-{}^{\rm t\!}Z
&O&
\end{pmatrix},
 \end{equation}
 where $Z=Z(t)=(z_1(t),\cdots,z_n(t))$, 
 $W=W(t)=(w_1(t),\cdots,w_n(t))$, 
 and for $F_{\frak m}$ in (4.40) with $F_{\frak k}\equiv0$, 
 the biharmonic map equation $(4.6)$ 
 is equivalent to 
 \begin{equation}
 \left\{
 \begin{aligned}
 -Z'''&-(\vert Z\vert^2+\vert W\vert^2)Z\\
 &+(
 2\langle Z,Z'\rangle+2\langle W,W'\rangle-\langle Z',Z\rangle
 -\langle W',W\rangle)Z\\
 &+(\langle Z',\overline{W}\rangle-\langle W',\overline{Z}\rangle)
 \overline{W}=0,\\
 -W'''&-(\vert Z\vert^2+\vert W\vert^2)W\\
 &+(
 2\langle Z,Z'\rangle+2\langle W,W'\rangle-\langle Z',Z\rangle
 -\langle W',W\rangle)W\\
 &+3(\langle Z',\overline{W}\rangle-\langle W',\overline{Z}\rangle)
 \overline{Z}=0,
 \end{aligned}
 \right.
 \end{equation}
 where 
 $Z'=(z_1{}'(t), \cdots, z_n{}'(t))$ and 
 $\langle Z,W\rangle:=\sum_{i=1}^nz_i(t)\,\overline{w_i(t)}$. 
 \par
 We find the following solutions of (4.41): 
 \par
 Case $(i)$: $z_1(t)=\cdots=z_n(t)=D_t $ and $w_1(t)=\cdots=w_n(t)=0$.
 \par
 Case $(ii)$: $z_1(t)=\cdots=z_n(t)=\sqrt{-1} 
 D_t $ and $w_1(t)=\cdots=w_n(t)=0$.
 \par
 Case $(iii)$: $z_1(t)=\cdots=z_n(t)=0$ and $w_1(t)=\cdots=w_n(t)=D_t$. 
 \par
 Case $(iv)$: $z_1(t)=\cdots=z_n(t)=0$ and 
 $w_1(t)=\cdots=w_n(t)=\sqrt{-1} D_t$. 
 \par
 The corresponding biharmonic curves into the quaternion projective spaces ${\mathbb H}P^n$ are given as follows: 
 \par
 Case $(i)$: 
 $$
 \varphi(t)=x\,\left[\cos (\sqrt{n}\,d_t),-\frac{1}{\sqrt{n}}\,\sin(\sqrt{n}\,d_t),\cdots,
 -\frac{1}{\sqrt{n}}\,\sin(\sqrt{n}\,d_t)\right].
 $$
 \par Case $(ii)$: 
 $$
 \varphi(t)=x\,\left[\cos (\sqrt{n}\,d_t),i\frac{1}{\sqrt{n}}\,\sin(\sqrt{n}\,d_t),\cdots,
 i\frac{1}{\sqrt{n}}\,\sin(\sqrt{n}\,d_t)\right].
 $$
 \par  Case $(iii)$: 
 $$
 \varphi(t)=x\,\left[\cos (\sqrt{n}\,d_t),-j\frac{1}{\sqrt{n}}\,\sin(\sqrt{n}\,d_t),\cdots,
 -j\frac{1}{\sqrt{n}}\,\sin(\sqrt{n}\,d_t)\right].
 $$
\par Case $(iv)$: 
 $$
 \varphi(t)=x\,\left[\cos (\sqrt{n}\,d_t),k\frac{1}{\sqrt{n}}\,\sin(\sqrt{n}\,d_t),\cdots,
 k\frac{1}{\sqrt{n}}\,\sin(\sqrt{n}\,d_t)\right].
 $$
Here, $x\in Sp(n+1)$, $i$, $j$ and $k$ are the quaternions 
satisfying $i^2=j^2=k^2=-1$ and $ij=k$, and 
$d_t=\frac{a}{3}t^3+\frac{b}{2}t^2+ct$, $a$, $b$ and $c$ are constant real numbers.   In each case, $\varphi$ is harmonic if and only if $a=b=0$. 
\vskip0.6cm\par
\section{Biharmonic maps from plane domains}
\subsection{Setting and deriving the equations}
In this section, we will treat with biharmonic maps of $(M,g)$ into a Riennian symmetric space $(G/K,h)$, with $\dim M=2$. We assume that 
$(M,g)=(\Omega,g)$ is 
an open domain in the 2-dimensional Euclidean space 
${\mathbb R}^2$ with $g=\mu^2\,g_0$, where
$\mu$ is a positive $C^{\infty}$ function on $\Omega$, 
$g_0=(dx)^2+(dy)^2$ is  the standard Euclidean metric and 
$(x,y)$  is the standard coordinate on ${\mathbb R}^2$. 
\par
Let $\varphi$ be a $C^{\infty}$ map from $\Omega$ into a symmetric space $N=G/K$ with a local lift $\psi:\,\Omega\rightarrow G$ satisfying 
$\varphi=\pi\circ\psi$, 
where $\pi:\,G\rightarrow G/K$ is the standard projection. 
The pull back of the Maurer-Cartan form $\theta$ on $G$ by $\psi$ 
is given by 
\begin{align*}
\alpha=\psi^{-1}d\psi &=\psi^{-1}\frac{\partial \psi}{\partial x}dx
+\psi^{-1}\frac{\partial \psi}{\partial y}dy\\
&=A_x\,dx+A_y\,dy,
\end{align*}
where we decompose two $\frak g$-valued functions
$A_x:=\psi^{-1}\frac{\partial \psi}{\partial x}$ and 
$A_y:=\psi^{-1}\frac{\partial \psi}{\partial y}$ 
 on $\Omega$ 
 according to the Cartan decomposition 
 ${\frak g}={\frak k}\oplus{\frak m}$ as follows: 
  $$
  A_x=A_{x,\,{\frak k}}+A_{x,\,{\frak m}},\quad 
  A_y=A_{y,\,{\frak k}}+A_{y,\,{\frak m}},
  $$
  which yield the decomposition of $\alpha$: 
  $
  \alpha=\alpha_{\frak k}+\alpha_{\frak m}, 
  $
  where 
  $$
  \alpha_{\frak k}=A_{x,\,{\frak k}}\,dx+A_{y,\,{\frak k}}\,dy, 
  \quad 
  \alpha_{\frak m}=A_{x,\,{\frak m}}\,dx+A_{y,\,{\frak m}}\,dy. 
  $$
  \par
  Then, we have by a direct computation, 
  \begin{equation}
  \delta(\alpha_{\frak m})=-\mu^{-2}\left\{
  \frac{\partial A_{x,\,{\frak m}}}{\partial x}
  + \frac{\partial A_{y,\,{\frak m}}}{\partial y}
  \right\}.
  \end{equation}
  Indeed, 
  if we take, as an orthonormal frame field with respect to $g$,  
  $e_1=\frac{1}{\mu}\frac{\partial}{\partial x}$ and 
  $e_2=\frac{1}{\mu}\frac{\partial}{\partial y}$. Then, we have 
  $$
  \alpha_{\frak m}(\nabla_{e_1}e_1)
  =-\mu^{-3}\frac{\partial\mu}{\partial y}A_{y,\,\frak m}, \quad 
  \alpha_{\frak m}(\nabla_{e_2}e_2)
  =-\mu^{-3}\frac{\partial\mu}{\partial x}A_{x,\,\frak m}, 
  $$ 
  and 
  $$
  \delta(\alpha_{\frak m})=
  -\sum_{i=1}^2\left\{
  \nabla_{e_i}(\alpha_{\frak m}(e_i))-\alpha_{\frak m}(\nabla_{e_i}e_i
  \right\},
  $$
  we have (5.1).\qed
  \vskip0.6cm\par
  Next, we have to calculate the harmonic map equation $(3.25)$, and the biharmonic map equation $(3.26)$ in this case. 
\par
First, 
for the left hand side of $(3.25)$, we have 
\begin{align}
&-\delta(\alpha_{\frak m})
+\sum_{i=1}^m [\alpha_{\frak k}(e_i),\alpha_{\frak m}(e_i)]\nonumber\\
&=\mu^{-2}
\left\{
\frac{\partial A_{x,{\frak m}}}{\partial x}+
\frac{\partial A_{y,{\frak m}}}{\partial y}
\right\}
+[\mu^{-1}\,A_{x,{\frak k}},\mu^{-1}\,A_{x,{\frak m}}]
+[\mu^{-1}\,A_{y,{\frak k}},\mu^{-1}\,A_{y,{\frak m}}]\nonumber\\
&=
\mu^{-2}
\left\{
\frac{\partial A_{x,{\frak m}}}{\partial x}+
\frac{\partial A_{y,{\frak m}}}{\partial y}
+[A_{x,{\frak k}},A_{x,{\frak m}}]
+[A_{y,{\frak k}},A_{y,{\frak m}}]
\right\}. 
\end{align}
For the left hand side of $(3.26)$, 
since 
$\Delta_g=-\mu^{-2}
\left\{
\frac{\partial^2}{\partial x^2}+\frac{\partial^2}{\partial y^2}
\right\}$, 
we have 
 \begin{align}
& \Delta_g(
 -\delta(\alpha_{\frak m})
+\sum_{i=1}^m [\alpha_{\frak k}(e_i),\alpha_{\frak m}(e_i)]
 )
 \nonumber\\
&\qquad +\sum_{s=1}^m
 [[
  -\delta(\alpha_{\frak m})
+\sum_{i=1}^m [\alpha_{\frak k}(e_i),\alpha_{\frak m}(e_i)],\alpha_{\frak m}(e_s)],\alpha_{\frak m}(e_s)]
 \nonumber\\
 &= -
 \mu^{-2}
\big\{
\frac{\partial^2}{\partial x^2}+\frac{\partial^2}{\partial y^2}
\big\}
\big(
\mu^{-2}
\big\{
\frac{\partial A_{x,{\frak m}}}{\partial x}+
\frac{\partial A_{y,{\frak m}}}{\partial y}
\nonumber\\
&
\qquad\qquad\qquad\qquad\qquad\quad\qquad+[A_{x,{\frak k}},A_{x,{\frak m}}]
+[A_{y,{\frak k}},A_{y,{\frak m}}]
\big\}
\big)\nonumber\\
&+
\mu^{-4}
\big[\big[
\frac{\partial A_{x,{\frak m}}}{\partial x}+
\frac{\partial A_{y,{\frak m}}}{\partial y}
+[A_{x,{\frak k}},A_{x,{\frak m}}]
+[A_{y,{\frak k}},A_{y,{\frak m}}], 
A_{x,{\frak m}}\big],A_{x,{\frak m}}\big]\nonumber\\
&+
\mu^{-4}
\big[\big[
\frac{\partial A_{x,{\frak m}}}{\partial x}+
\frac{\partial A_{y,{\frak m}}}{\partial y}
+[A_{x,{\frak k}},A_{x,{\frak m}}]
+[A_{y,{\frak k}},A_{y,{\frak m}}], 
A_{y,{\frak m}}\big],A_{y,{\frak m}}\big]. 
 \end{align}
 Therefore, we have that $\varphi:\,(\Omega,g)\rightarrow (G/K,h)$ is biharmonic if and only if the right hand side of $(5.3)$ vanish. 
 \par
 Second, we have to see the integrability condition $(3.23)$ or $(3.24)$. We have 
 \begin{align}
 d\alpha_{\frak k}&+\frac12[\alpha_{\frak k}\wedge\alpha_{\frak k}]+\frac12[\alpha_{\frak m}\wedge\alpha_{\frak m}]\nonumber\\
 &=
 \left\{
 -\frac{\partial A_{x,\,{\frak k}}}{\partial y}
 +\frac{\partial A_{y,\,{\frak k}}}{\partial x}
 +[A_{x,\,{\frak k}}, A_{y,\,{\frak k}}]
 +[A_{x,\,{\frak m}}, A_{y,\,{\frak m}}]
 \right\}
 \,dx\wedge dy\nonumber\\
 &=0,\nonumber
 \end{align}
 so that we have 
 \begin{equation}
  -\frac{\partial A_{x,\,{\frak k}}}{\partial y}
 +\frac{\partial A_{y,\,{\frak k}}}{\partial x}
 +[A_{x,\,{\frak k}}, A_{y,\,{\frak k}}]
 +[A_{x,\,{\frak m}}, A_{y,\,{\frak m}}]=0.
 \end{equation}
 For the second equation of $(3.24)$, 
 $d\alpha_{\frak m}+[\alpha_{\frak k}\wedge \alpha_{\frak m}]=0$, 
 we have 
  \begin{equation}
  -\frac{\partial A_{x,\,{\frak m}}}{\partial y}
 +\frac{\partial A_{y,\,{\frak m}}}{\partial x}
 +[A_{x,\,{\frak k}}, A_{y,\,{\frak m}}]
 +[A_{x,\,{\frak m}}, A_{y,\,{\frak k}}]=0.
 \end{equation}
 Summing up the above, we obtain 
 \begin{th}
 Let $\Omega\subset {\mathbb R}^2$ an open domian, $g=\mu^2g_0$, 
 $\mu>0$, a positive $C^{\infty}$ function on $\Omega$, and 
 $g_0=(dx)^2+(dy)^2$ is the standard Riemannian metric 
 on ${\mathbb R}^2$. 
 on which $(x,y)$ is the standard coordinate. 
 Let  $(G/K,h)$ a Riemannian symmetric space, 
 with $\pi:\,G\rightarrow G/K$, the projection. 
 For every $C^{\infty}$ map from $\Omega$ into $G/K$ with a local lift 
 $\psi:\,\Omega\rightarrow G$ such that 
 $\varphi=\pi\,\circ\,\psi$, 
 let $\alpha=\psi^{\ast}\theta$, the pull back of the Maurer-Cartan form 
 $\theta$ on $G$ by $\psi$ and decompose it in such a way that 
 $\alpha=\alpha_{\frak k}+\alpha_{\frak m}$ corresponding to the Cartan decomposition $\frak g=\frak k\oplus\frak m$. Then, 
 \par $(1)$ 
 $\varphi:\,(\Omega,g)\rightarrow (G/K,h)$ is harmonic if and only if 
 \begin{equation}
\frac{\partial A_{x,{\frak m}}}{\partial x}+
\frac{\partial A_{y,{\frak m}}}{\partial y}
+[A_{x,{\frak k}},A_{x,{\frak m}}]
+[A_{y,{\frak k}},A_{y,{\frak m}}]=0.
\end{equation}
\par $(2)$ 
 $\varphi:\,(\Omega,g)\rightarrow (G/K,h)$ is biharmonic if and only if 
 $(5.3)$ vanishes. 
 \par $(3)$ 
 For the integrability condition, $(5.4)$ and $(5.5)$ must hold. 
 \par
$(4)$  In particular, for a {\em horizontal} lift $\psi$, i.e., $\alpha_{\frak k}\equiv 0$,  we have 
 \begin{align}
 &-\bigg\{
 \frac{\partial^2}{\partial x^2}+ \frac{\partial^2}{\partial y^2}
 \bigg\}
 \bigg(\mu^{-2}\bigg\{
 \frac{\partial P}{\partial x}+\frac{\partial Q}{\partial y}
 \bigg\}\bigg)
 +\bigg[\bigg[
 \mu^{-2}\bigg\{
 \frac{\partial P}{\partial x}+\frac{\partial Q}{\partial y}\bigg\},P\bigg],P\bigg]\nonumber\\
 &\qquad\quad\qquad\quad+\bigg[\bigg[
 \mu^{-2}\bigg\{
 \frac{\partial P}{\partial x}+\frac{\partial Q}{\partial y}\bigg\},Q\bigg],Q\bigg]=0,\\
& [P,Q]=0,\\
&-\frac{\partial P}{\partial y}+\frac{\partial Q}{\partial x}=0,
 \end{align}
 where we put 
 $P:=\alpha_{x,\,{\frak m}}$ and $Q:=\alpha_{y,\,{\frak m}}$. 
 In the case $\mu=1$, the following three equations must hold for the biharmonic map $\varphi$: 
 \begin{align}
 &-P_{xxx}-P_{xyy}-Q_{xxy}-Q_{yyy}\nonumber\\
 &\qquad+[[P_x+Q_y,P],P]+[[P_x+Q_y,Q],Q]=0,\\
 &[P,Q]=0,\\
 &P_y-Q_x=0,
 \end{align}
 where we denote $P_x=\frac{\partial P}{\partial x}$, etc. 
 \end{th}
 \vskip0.6cm\par
 \subsection{Solving the biharmonic map equations} 
 In this subsection, we want to give the solutions of the equations $(5.10)$, 
 $(5.11)$ and $(5.12)$. 
 \par
 To do it, let us consider the special case that 
 $P_y\equiv0$ and $Q_x\equiv0$, i.e., 
 $P(x,y)=P(x)$ and $Q(x,y)=Q(y)$.  Then, $(5.12)$ holds clearly. 
 The left hand side of $(5.10)$ coincides with 
 \begin{align}
& \big\{
 -P_{xxx}+\big[\big[P_x,P\big],P\big]
 \big\}
 +\big\{
 -Q_{yyy}+\big[\big[Q_y,Q\big],Q\big]
 \big\}\nonumber\\
 &+\big[\big[Q_y,P\big],P\big]+\big[\big[P_x,Q\big],Q\big]=0.
 \end{align}
 Here, we have that 
 $[[Q_y,P],P]=0$ and $[[P_x,Q],Q]=0$. 
 Because, we have due to $Q_x=0$
 \begin{equation}
 \frac{\partial}{\partial x}[[P,Q],Q]
 =[[P_x,Q],Q].
 \end{equation}
 But, due to $(5.11)$ the left hand side of $(5.14)$ must vanish. 
 By the same way, we have $[[Q_y,P],P]=0$. 
 \par
 Thus, $(5.13)$ turns out that 
 \begin{equation}
  -\big\{-P_{xxx}+\big[\big[P_x,P\big],P\big]\big\}
  =-Q_{yyy}+\big[\big[Q_y,Q\big],Q\big]
 \end{equation}
 But, notice that 
 the left hand side of $(5.15)$ is an $\frak m$-valued function only in $x$ and 
 the right hand side of $(5.15)$ is the one only in $y$, so we have 
 \begin{equation}
 \left\{\begin{aligned}
 -P_{xxx}+\big[\big[P_x,P\big],P\big]&=c,\\
 -Q_{yyy}+\big[\big[Q_y,Q\big],Q\big]&=-c,
 \end{aligned}
 \right.
 \end{equation}
 where $c\in {\frak m}$ is a constant vector. 
 \par
 Notice here that both two equations of $(5.16)$ are the same as $(4.6)$ in the case $c=0$. So, we can obtain the following two theorems by carrying out the similar calculations as in 4.2. 
 \par
 Thus, we have
 \begin{th}
 Let $(G/K,h)$ be a Riemannian symmetric space whose rank 
 is bigger than or equal to two, $\frak g=\frak k\oplus\frak m$, 
 the Cartan decomposition, 
 $\frak a$, a maximal abelian subalgebra of $\frak g$ 
 contained in $\frak m$. 
 Let $X,\,Y\in {\frak a}$  be two elements in $\frak a$ which are linearly independent. 
 \par
 $(1)$ Let us take two $\frak m$-valued functions $P(x,y)=(a_1\,x^2+b_1\,x+c_1)\,X$ and 
 $Q(x,y)=(a_2\,y^2+b_2\,y+c_2)\,Y$, where 
 $a_i$, $b_i$ and $c_i$ $(i=1,2)$ are constant real numbers. 
 Then, $P$ and $Q$ are solutions of $(5.10)$, $(5.11)$ and $(5.12)$. 
 For such $P$ and $Q$, there exists a unique $C^{\infty}$ map $\psi$ 
 from $\Omega$ into $G$ such that 
 $\varphi=\pi\,\circ\,\psi$ is a biharmonic mapping 
 form $(\Omega,g_0)$ into $(G/K,h)$ with $\varphi(0,0)=x_0\in G$ 
 for a fixed point $x_0\in G/K$. 
 $\varphi:\,(\Omega,g)\rightarrow (G/K,h)$ is harmonic if and only if 
 $a_i=b_i=0$ $(i=1,2)$. 
 \par
 $(2)$ Assume that $G$ is a matrix Lie group, i.e., a subgroup of $GL(N,{\mathbb C})$. Then, the above $C^{\infty}$ maps $\psi:\,\Omega\rightarrow G$ and 
 $\varphi=\pi\,\circ\,\psi$ are given by 
 \begin{equation}
 \left\{
 \begin{aligned}
 \psi(x,y)&=x_0\,\exp(d_x\,X+d_y\,Y)\in G, 
 \\
 \varphi(x,y)&=x_0\,\exp(d_x\,X+d_y\,Y)\,\cdot\,o\in G/K,
 \end{aligned}
 \right.
 \end{equation}
where $o=\{K\}\in G/K$, $d_x=\frac{a_1}{3}\,x^3+\frac{b_1}{2}\,x^2+c_1\,x$ and 
$d_y=\frac{a_2}{3}\,y^3+\frac{b_2}{2}\,y^2+c_2\,y$, respectively. 
 \end{th}
 \begin{pf}
 We only have to see $(2)$. 
 By the assumption that $\{X,Y\}$ is abelian, we have 
 for the $\psi(x,y)$ of the form $(5.17)$, as a matrix of degree $N$, 
 \begin{align*}
 \frac{\partial \psi}{\partial x}&=x_0\,\exp(d_x\,X+d_y\,Y)\,\cdot\,
 \frac{\partial}{\partial x}(d_x\,X+d_y\,Y)\\
 &=\psi\,\cdot\,(a_1\,x^2+b_1\,x+c_1)\,X\\
 &=\psi\,P,
 \end{align*}
 so we have $\psi^{-1}\frac{\partial \psi}{\partial x}=P$. By the same way,  $\psi^{-1}\frac{\partial \psi}{\partial y}=Q$, so 
 we have 
 $\psi^{-1}d\psi=P\,dx+Q\,dy=\alpha$. The mapping $\psi$ is the desired 
 $C^{\infty}$ mapping of $\Omega$ into $G$, and due to Theorem 3.6, 
 we obtain a biharmonic mapping of $(\Omega, g_0)$ into $(G/K,h)$. 
 \end{pf}
 \vskip0.6cm\par
 {\it Remark.} \quad When $a_i=b_i=0$ $(i=1,2)$, 
 the mapping $\varphi:\,{\mathbb R}^2\rightarrow(G/K,h)$ is a 
 well known totally geodesic immersion into a Riemannian symmetric space $(G/K,h)$. 
 \vskip0.6cm\par
 By the similar calculation as the subsection 4.2, we obtain 
 \begin{th}
 For the cases of the standard unit sphere $(S^n,h)$, the complex projective space 
 $({\mathbb C}P^n,h)$, the quaternion one $({\mathbb H}P^n,h)$, we obtain the following biharmonic mappings 
 of $({\mathbb R}^2,g_0)$ into them, respectively. 
 \par
 $(1)$ Case of $(S^n,h)$: 
 \begin{align}
  \varphi_1(t)
 =&x_0\,{}^{\rm t\!}(\cos(\sqrt{n} (d_x+d_y)),\nonumber\\
& \frac{1}{\sqrt{n}}\sin(\sqrt{n}(d_x+d_y)),\cdots,
\frac{1}{\sqrt{n}}\sin(\sqrt{n}( d_x+d_y))),
 \end{align}
is a biharmonic mapping of $({\mathbb R}^2,g_0)$ into $(S^n,h)$, where $x_0\in G=SO(n+1)$. 
\par
$(2)$ Case of $({\mathbb C}P^n,h)$: 
 \begin{align}
  \varphi_2(t)
 =&x_0\,{}^{\rm t\!}(\cos(\sqrt{n} (d_x+d_y)),\nonumber\\
& \frac{\sqrt{-1}}{\sqrt{n}}\sin(\sqrt{n}(d_x+d_y)),\cdots,
\frac{\sqrt{-1}}{\sqrt{n}}\sin(\sqrt{n}( d_x+d_y))),
 \end{align}
is a biharmonic mapping of 
$({\mathbb R}^2,g_0)$ into $({\mathbb C}P^n,h)$, 
where $x_0\in G=SU(n+1)$. 
\par
$(3)$ Case of $({\mathbb H}P^n,h)$: 
 \begin{align}
  \varphi_3(t)
 =&x_0\,{}^{\rm t\!}(\cos(\sqrt{n} (d_x+d_y)),\nonumber\\
& \frac{k}{\sqrt{n}}\sin(\sqrt{n}(d_x+d_y)),\cdots,
\frac{k}{\sqrt{n}}\sin(\sqrt{n}( d_x+d_y))),
 \end{align}
is a biharmonic mapping of 
$({\mathbb R}^2,g_0)$ into $({\mathbb H}P^n,h)$, 
where $x_0\in G=Sp(n+1)$, and 
$i$, $j$ and $k$ are the quaternions satisfying 
$i^2=j^2=k^2=-1$ and $ij=k$. 
\par
Here, in all the cases, 
$d_t=\frac{a}{3}t^3+\frac{b}{2}t^2+ct$ for $t=x$ or $t=y$.   
\par
Furthemore, each $\varphi_i$ $(i=1,2,3)$ are harmonic of $({\mathbb R}^2,g_0)$ into 
$(S^n,h)$, $({\mathbb C}P^n,h)$ or $({\mathbb H}P^n,h)$ if and only if $a=b=0$, respectively. 
 \end{th}
\vskip2cm\par       

\end{document}